\documentclass[11pt]{article}
\usepackage{anysize}
\marginsize{1in}{1in}{1in}{1in}
\usepackage{amsfonts,amsmath,amsxtra,euscript}
\usepackage{amssymb}
\usepackage{graphicx,longtable,wrapfig}
\usepackage{epstopdf}

\usepackage[english]{babel}

\author{Alexander Esterov\thanks{National Research University Higher School of Economics \newline This study was carried out within The National Research University Higher School of Economics Academic Fund Program in 2012-2013, research grant No.11-01-0125. Partially supported by RFBR
grants 10-01-00678 and 12-01-31233, MESRF grant MK-6223.2012.1, and the Dynasty Foundation fellowship.}, Gleb Gusev\thanks{Moscow Institute of Physics and Technology (State University)\newline Partially supported by RFBR
grants 10-01-00678, 12-01-31233 and 13-01-00755.}}
\title{Systems of equations with a single solution}
\date{}

\newtheorem{thm}{Theorem}
\newtheorem{cor}{Corollary}
\newtheorem{prop}{Proposition}
\newtheorem{lem}{Lemma}

\newenvironment{dfntn} {\smallskip\noindent{\bf
Definition\/}.}{\smallskip\par}

\newenvironment{exmpl} {\smallskip\noindent{\bf Example\/}.}{\smallskip\par}
\newenvironment{rmrk} {\smallskip\noindent{\bf Remark\/}.}{\smallskip\par}
\newenvironment{prf} {\noindent{\em Proof\/}.}{{ $\Box$}\smallskip\par}

\renewcommand{\a}{\alpha }
\newcommand{\Ps}{\Psi }

\renewcommand{\d}{\delta }
\newcommand{\D}{\Delta }

\renewcommand{\l}{\lambda }

\newcommand{\OOO}{\Omega}
\newcommand{\MV}{{\rm MV}}

\newcommand{\R}{\mathbb{R}}

\newcommand{\C}{\mathbb{C}}

\newcommand{\Z}{\mathbb{Z}}
\newcommand{\ZZ}{\EuScript{Z}}
\renewcommand{\AA}{\widetilde{A}}

\newcommand{\Vol}{\mathop{\mathrm{Vol}}\nolimits}

\begin{document}
\maketitle

\begin{abstract}
We classify generic systems of polynomial equations with a single solution, or, equivalently, collections of lattice polytopes of minimal positive mixed volume. As a byproduct, this classification provides an algorithm to evaluate
the single solution of such a system.
\end{abstract}

\section{Introduction}

The mixed volume is the unique symmetric real-valued function $\MV$ of $n$ convex bodies in an $n$-dimensional vector space $V$
such that
$$
\MV(A+B,A_2,\ldots,A_n)=\MV(A,A_2,\ldots,A_n)+\MV(B,A_2,\ldots,A_n)
$$
in the sense of Minkowski summation $A+B=\{a+b\,|\, a\in A,\, b\in B\}$,
and $\MV(A,\ldots,A)=($volume of $A)$ in the sense of a given volume form on $V$. In what follows, $V$ is always of the form $\mathcal{V}\otimes\R$, where $\mathcal{V}$ is an integer lattice, and the volume form is always chosen in such a way 
that the volume of the torus $V/\mathcal{V}$ is equal to $n!$. For this volume form, the mixed volume of lattice polytopes (i.e. 
the ones whose vertices are in $\mathcal{V}$) is always integer.

The mixed volume is always non-negative. A criterion for a collection of polytopes to have zero mixed volume was established by Minkowski in 1911 in the original paper \cite{mink}, where multidimensional mixed volumes were introduced: the mixed volume of bodies is positive if and only if they contain linearly independent segments. A constructive version of this criterion was given by D.~Bernstein and A.~Khovanskii (\cite{kh77}, see ,e.g.,
\cite[Lemma 1.2]{e10} 
for a proof):
\begin{prop}\label{prop}
The mixed volume of convex bodies is zero if and only if $k$ of these bodies sum up to a body of dimension strictly smaller than $k$.
\end{prop}

In this paper, we provide the classification of collections of lattice polytopes with the minimal positive mixed volume (by induction on $n$):


\begin{thm} A collection $A$ of $n$ lattice polytopes in $V$ has the unit mixed volume if and only if

\noindent 1) the mixed volume is not zero, and

\noindent 2) there exists $k>0$
such that, up to translations, 
$k$ of the polytopes are faces of the same $k$-dimensional volume 1 lattice simplex in a $k$-dimensional rational subspace $U\subset V$, and the images of the other $n-k$ polytopes under the projection $V\to V/U$ 
have the unit mixed volume.
\end{thm}
See Section~\ref{Theorem-1-reduction} for the proof. Here and in what follows, the volume forms in the subspace $U\subset V$ and in the quotient space $V/U$ are induced by the lattices $\mathcal{U}=U\cap\mathcal{V}$ and $\mathcal{V}/\mathcal{U}$ respectively. Condition (1) implies that the $k$ polytopes in (2) generate the whole simplex.

\begin{exmpl} Any pair of lattice polygons of mixed area 1 is equal (up to translations and an authomorphism of $\Z^2$) to exactly one of the following pairs (with $a\geqslant b\geqslant 0$):

\begin{center}
\noindent\includegraphics[width=9.5cm]{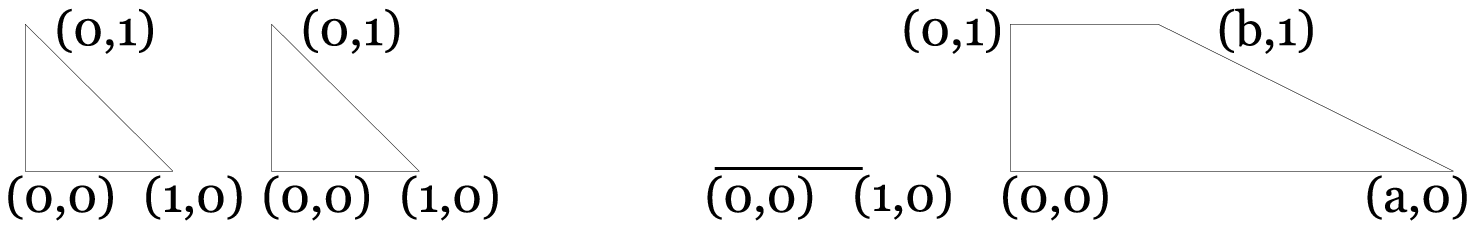}
\end{center}

\end{exmpl}

The question of classifying lattice polytopes of small mixed volume is particularly motivated by the study of codimensions of discriminants (see, e.g., \cite[Theorem 3.13]{me8}, 
or \cite{sturmf11}  for details). Theorem 1 was conjectured in \cite[Conjecture 3.16]{me8}, 
and its special case of full-dimensional polytopes was proved in \cite[Proposition 2.7]{sturmf11}.

The mixed volume is related to algebra by the Kouchnirenko-Bernstein formula: a system of $n$ polynomial equations of $n$ variables with Newton polytopes $N_1,N_2\ldots,N_n$ 
and generic
coefficients has $\MV(N_1,N_2,\ldots,N_n)$ 
solutions in the complex torus $(\C\setminus 0)^n$, see \cite{bernst}. Thus, Theorem 1 classifies all generic
systems of polynomial equations with a unique solution. By a general Gr\"{o}bner basis or Galois theory argument, the solution of such a system admits a rational expression in terms of the coefficients of the system, and Theorem 1 provides an explicit construction for it (by induction on $n$):

I) upon a certain monomial change of variables, 
$k$ of the $n$ equations become linear (non-homogeneous)
equations of $k$ variables, from which the $k$ variables can be evaluated;

II) after the substitution of the evaluated variables in the other equations, we obtain $n-k$ generic
equations of $n-k$ variables with a unique solution and 
proceed to the next group of simultaneously linearizable equations.

\begin{exmpl} In order to solve a system $a+bxy=0,\, f(xy)+yg(xy)=0$, we first make a monomial change $xy=u,\, y=v$, then solve the first equation $a+bu=0$, linear in $u$, and obtain $u=-a/b$, then put the result to the second equation $f(-a/b)+vg(-a/b)=0$, linear in $v$, and obtain $v=-f(-a/b)/g(-a/b)$. Theorem 1 ensures that such an obvious approach works in any dimension for square systems with a single solution. \end{exmpl}

This algorithm addresses the simplest possible problem in the field of efficient solving of polynomial systems with few monomials/solutions (see, e.g., \cite{few1} and \cite{few2} for motivation and recent advances in some other related problems). A polynomial-time implementation of this scheme is described in Section \ref{Sshift}. An example of further classification (square systems of $\leq$3 equations with two solutions) is given in Section \ref{Smv2}.


\section{Reduction of Theorem 1 to essential tuples}\label{Theorem-1-reduction}

\begin{dfntn}
A set of lattice
polytopes
$A_1, A_2,\ldots, A_k\subset \R^n$ is called {\it a $k$-tuple in $\R^n$}. Depending on the context, we treat a tuple $A$ either as an unordered multiset (when the order of its elements does not matter) or as a map $A\colon \{1,2,\ldots, k\}\to \OOO_n$ to the set $\OOO_n$ of lattice polytopes in $\R^n$ such that $A(j)=A_j$. For a tuple $A$, its subset $B\subset A$ is called {\it a subtuple} of $A$. Formally, if $A\colon \{1,2,\ldots, k\}\to \OOO_n$ is a map, its subtuple $B$ is the restriction of $A$ to a subset of $\{1,2,\ldots, k\}$.
\end{dfntn}

\begin{dfntn}
A $k$-tuple $A$ in $\R^n$ is
said to be {\it essential}, if for each $l\leq k$ and $j_1<j_2<\ldots<j_l$,
we have $\dim \sum_{s=1}^l A(j_s) \geq \min(l+1,n)$.
\end{dfntn}

For instance, any tuple in $\R^n$ that contains a segment is not essential, if $n>1$. See \cite{Sturmfels} for a relation between this notion and codimension of resultants. We now reformulate Theorem 1 as follows.

\begin{thm}\label{Thm1}
The mixed volume of an essential $n$-tuple $A$ in $\R^n$ is 1, if and only if all the polytopes
$A(i)$, up to translations,
are contained in a volume 1 lattice simplex. These translations are unique up to a simultaneous translation.
\end{thm}

The proof is given in Sections~\ref{preliminary}, \ref{exists} (existence of the translations) and Section~\ref{unique} (uniqueness).

\vspace{2ex}

\noindent{\it Proof of Theorem 1.} For an essential tuple $A$, Theorem 1 immediately follows from Theorem \ref{Thm1}. If the tuple $A$ is not essential, then it contains a minimal subtuple of $k$ polytopes that sums up to a $k$-dimensional polytope. In turn, this subtuple is essential. We apply Theorem~\ref{Thm1} to this subtuple and deduce Theorem 1 from the following well known fact. $\quad\Box$

\begin{prop}\label{prop-2} If the polytopes $A_1,A_2\ldots,A_k$ of an $n$-tuple $A=\{A_1,A_2,\ldots, A_n\}$ in $\R^n$ are contained in a subspace $\R^k\subset\R^n$, then
$$
\MV(A_1,A_2\ldots,A_n)=\MV(A_1,A_2\ldots,A_k)\cdot\MV(B_{k+1},B_{k+2}\ldots,B_n),
$$
where $B_j$ is the image of $A_j$ under the projection along $\R^k$.
\end{prop}

Note that the uniqueness statement in Theorem \ref{Thm1} fails for non-essential tuples of subsimplices of a volume 1 lattice simplex; the worst case in this respect are the unit segments $e_1,\ldots,e_n$ on the coordinate axes in $\Z^n$:

\begin{prop} \label{cayley} There exist exactly $2^n(n+1)^{n-2}$ lattice volume 1 simplices, containing the segments $e_1,\ldots,e_n$ up to a translation.\end{prop}

The proof is similar to the Cayley tree-counting formula, because the desired simplices are in one to one correspondence with the trees that consist of translated copies of $e_1,\ldots,e_n$. Note that this quantity also equals the number of $n$-trees in the $n$-dimensional cube.

\vspace{2ex}

The key ideas in the proof of Theorem \ref{Thm1} are the use of the Aleksandrov--Fenchel inequality
\begin{equation}\label{AF}
\MV(A,B,C,\ldots)^2\geqslant\MV(A,A,C,\ldots)\cdot\MV(B,B,C,\ldots)
\end{equation}
and the notion of %
the mixed fiber body:

\begin{dfntn} The mixed fiber body of the convex bodies $A_1,A_2\ldots,A_{n-k}\subset V$ 
in a $(k+1)$-dimensional rational subspace $U\subset V$ is the convex body $X\subset U$ 
such that $\MV(X,B_1,\ldots,B_k)=\MV(A_1,\ldots,A_{n-k},B_1,\ldots,B_k)$ for all collections of convex bodies $B_1,B_2\ldots,B_k\subset U$.
\end{dfntn}
The mixed fiber body always exists, is unique up to a translation, and is a lattice polytope if so are $A_1,A_2\ldots,A_{n-k}$ 
(see~\cite{E-Kh} for the proofs and the equivalence to the original definition~\cite{mcd}).
For instance, if $A_1=\ldots=A_{n-k}=A$ is a volume 1 simplex whose $(k+1)$-dimensional face $F$ is parallel to $U$, then the mixed fiber body equals $F$ up to a translation.

Note that Theorem \ref{Thm1} can be extended to mixed fiber bodies as follows.
Consider a $k$-tuple $A$ in $\R^n$ and a rational affine subspace $L\subset \R^n$ of dimension $n-k+1$.

\begin{dfntn}
We call the tuple $A$ {\it essential with respect to $L$}, if, for each $(n-k)$-tuple $B$ of full-dimensional polytopes
in $L$, the tuple $A\cup B$ is essential.
\end{dfntn}

\begin{cor}\label{Cor1}
Assume that
a $k$-tuple $A$ in $\R^n$ with $0<k<n$ is essential with respect to a subspace $L\subset \R^n$ of dimension $n-k+1$, and the mixed fiber body $\AA \subset L$ of
the tuple $A$ is contained, up to translations,
in a volume 1 lattice full-dimensional simplex $\delta$ in $L$. Then there exists a volume 1 lattice
simplex $\Delta \subset \R^n$ that is adjacent to $\delta$ and contains all the polytopes $A(i)$ up to translations.
\end{cor}

\begin{prf}
Consider the tuple $B$ consisting of $n-k$ copies of
$\d$. The tuple $A\cup B$ is essential, and the mixed volume of $A\cup B$ equals the mixed volume of $\{\d\}\cup B$, which equals $1$. We apply Theorem~\ref{Thm1} to $A\cup B$ and obtain a required volume 1 simplex $\D$.
\end{prf}

Note that the condition $k<n$ in Corollary~\ref{Cor1} is essential. In fact, in the case $k=n$, any essential tuple $A$ is essential with respect to any one-dimensional subspace $L\subset \R^n$.
Let $A(i)=\Delta$ be a
volume~1 lattice full-dimensional simplex in $\R^n$ for $i=1,2,\ldots, n$.
Given an arbitrary $L$, the mixed fiber body $\AA\subset L$ is a unit segment. However, the simplex $\Delta$ is not adjacent to this segment in general.

\section{Preliminary statements and notations}\label{preliminary}

In what follows, we use the following notations throughout the paper:
\begin{itemize}
\item For a lattice polytope $X\in \OOO_n$, we define its support function
$l_X\colon (\Z^n)^* \to \Z$ by $l_X(\a) = \max_{x\in X}\a(x)$, where $\a\in (\Z^n)^*$.
\item  For an $l$-dimensional lattice polytope $X$, we use $\Vol(X)$ for its $l$-dimensional volume in the affine span of $X$.
\item We denote by $\ZZ\subset (\Z^n)^*$ the set of primitive integer covectors.
\item For a covector $\a\in \ZZ$, we denote by $X^\a$ the face $\{x\in X\mid \a(x)=l_X(\a)\}$ of the
polytope~$X$, where the covector $\a$ attains its maximal value.
\item For a $k$-tuple $A$ in a rational $l$-dimensional
affine subspace $L\subset \R^n$,

\begin{itemize}
 \item the Minkowski sum of the polytopes of tuple $A$ is denoted by $\sum A = \sum_{i=1}^k A(i)$;

\item if $k=l$, the mixed volume of $A$ is denoted by
$$
\prod A =\prod_{j=1}^k A(j)=\MV(A(1),A(2),\ldots, A(k));
$$

\item we define {\it dimension} of tuple $A$ as $\dim A:= \dim \sum A - k$ .

\end{itemize}
\end{itemize}

\begin{dfntn}
A $k$-tuple $A$ in $\R^n$ is said to be {\it linearly independent} if, for each $l$ and $j_1<j_2<\ldots<j_l$,
one has $\dim \sum_{s=1}^l A(j_s) \geq l$.
\end{dfntn}

The following identity (see \cite{bernst}) often turns out useful for treatment and computation of the mixed volume of an $n$-tuple $A$ in $\R^n$:

\begin{equation}\label{IV_reduce}
\prod A = \sum_{\a\in \ZZ} l_{A(1)}(\a) \prod_{j=2}^n A(j)^{\alpha}.
\end{equation}

As a consequence of this equation, the following {\it monotonicity property} of the mixed volume holds. For a couple of $n$-tuples $A,B$ in $\R^n$, if $A(j)\subset B(j)$, $j=1,2,\ldots, n$, then $\prod A\leq \prod B$.

The following two statements support the key inductive step in proof of Theorem~\ref{Thm1}.

\begin{lem}\label{Delta}

Let an $(n-1)$-tuple $A=\{A_1,A_2,\ldots, A_{n-1}\}$ in $\R^n$ be essential and consist
of lattice polytopes contained (up to translations) in a volume 1 lattice full-dimensional simplex $\Delta$ in $\R^n$.
\begin{enumerate}
\item For each covector $\a\in \ZZ$ such that $\D^{\a}$ is a facet of $\D$, the mixed volume $\prod (A_i)^{\a}$ equals~$1$.

\item For a covector $\a\in \ZZ$ such that $\dim \D^{\a} < n-1$, we have $\prod (A_i)^{\a}=0$.
\end{enumerate}

\end{lem}

\begin{prf}
1. Assume for simplicity that $A_i\subset \D$ and $(A_i)^{\a}\subset \D^{\a}$ for {$i=1,2,\ldots, n-1$}, which is actually true up to translations. Since the polytopes $B_i: = (A_i)^{\a}$ lie in the volume~1 simplex $\D^{\a}$, $i=1,2,\ldots, n-1$, we have $\prod B \leq \Vol (\D^{\a}) = 1$. It remains to prove that $\prod B \neq 0$.

Assume the contrary: $\prod B=0$. Proposition~\ref{prop} provides a subtuple
$$
B' = \{B_{j_1}, B_{j_2}, \ldots, B_{j_k}\},\, j_1<j_2<\ldots<j_k,
$$
of $B$ such that ${\dim \sum B' < k}$. Consider a graph that consists of $k$ nodes corresponding to the elements of $B'$ and has an edge between a pair of polytopes $B_{x_1}, B_{x_2}\in B'$ if and only if their intersection is not empty: $B_{x_1}\cap B_{x_2}\neq \emptyset$. Let $C_1, C_2,\ldots ,C_l$ be the subtuples of $B'$ corresponding to the connected components of the graph. Let $X_i$ be the number of vertices of $\D$ that belong to the union $\cup_{X\in C_j} X$ of polytopes
in $C_i$. We have $\sum_j X_j = n$. One can easily see that $\dim \sum B' = \sum_{j=1}^t (\dim \sum C_j)$ and $\dim \sum C_j = X_j-1$. Denote $c_t$ the number of polytopes in $C_t$. We have $\sum_{t=1}^l c_t=k> \dim \sum B'= \sum (X_i - 1)$. Thus there exists $t$ such that $X_t - 1 < c_t$. Let $B_{s_1}, B_{s_2}, \ldots, B_{s_{c_t}}$ be
the polytopes of $C_t$. Then
$$
\dim \sum_{j=1}^{c_t} A_{s_j} \leq 1+\dim \sum_{j=1}^{c_t} B_{s_j} = 1+\dim \sum C_t \leq X_t \leq c_t,
$$
which contradicts to the statement of Lemma~\ref{Delta} requiring $A$ to be essential.

2. Denote $\a_1,\a_2,\ldots \a_{n+1}\in \ZZ$ the covectors such that $\D^{\a_j}$ are the $n+1$ facets of $\D$. Assume there exists a covector $\a_0\in \ZZ$ different from $\a_j$ such that $\prod (A_j)^{\a_0}>0$. Consider any edge $(v_0,v_1)$ of $\D$ that is not orthogonal to $\a_0$, assume $\a_0(v_0)< \a_0(v_1)$. Choose an integer affine coordinate system in $\Z^n$ with the origin at $v_0$. W.l.o.g., assume that $\D^{\a_{n+1}}$ is the facet of $\D$ that does not contain $v_0$. By the identity~(\ref{IV_reduce}), one has $\prod (A\cup \{\D\}) = \sum_{\a\in \ZZ} l_{\Delta}(\a) \prod (A_j)^{\a} \geq l_{\Delta}(\a_{n+1}) \prod (A_j)^{\a_{n+1}} + l_{\Delta}(\a_0) \prod (A_j)^{\a_0} \geq 2$. On the other hand, the monotonicity property of mixed volumes implies that $\prod (A\cup \{\D\}) \leq \Vol (\D) = 1$, which leads to a contradiction.
\end{prf}

\begin{lem}\label{step_lem}
Let an $(n-1)$-tuple $A$ in $\R^n$ be
essential and consist
of lattice polytopes contained (up to translations)
in a volume 1 lattice full-dimensional
simplex $\Delta$
in $\R^n$. If the mixed volume $\prod (A\cup \{X\})$ with
some lattice polytope
$X$ equals 1, then $X$ is also contained in $\D$ up
to a translation.
\end{lem}

\begin{prf}
Chose a lattice basis such that $\D$ is the convex hull of the origin and the basis vectors. Let $\a_1,\a_2,\ldots, \a_{n+1} \in \ZZ$ be the covectors such that $\D^{\a_j}$ are the $n+1$ facets of $\D$. Assume that $\D^{\a_{n+1}}$ is the facet that does
not contain the origin. Let $v\in R^n$ be the vector defined by the conditions $\a_i(v) = l_{X}(\a_{i}), i=1,2\ldots, n$. Then, for the polytope
$X' = X+ \{-v\}$, one has
$l_{X'}(\a_{i}) = 0 = l_{\D}(\a_{i})$. Due to Lemma~\ref{Delta}, identity~(\ref{IV_reduce}) implies that $\prod (A\cup \{X'\}) = \sum_i l_{X'}(\a_{i})\prod_{j=1}^{n-1} (A_j)^{\a_i} = l_{X'}(\a_{n+1})$. On the other hand, $\prod (A\cup \{X'\}) = \prod (A\cup \{X\}) = 1$. Therefore, $l_{X'}(\a_{n+1}) = 1$, and $X' \subset \Delta$.
\end{prf}

This and the subsequent lemmas enable to reduce, under appropriate conditions, the prove of Theorem~\ref{Thm1} concerning a tuple $A$ to the prove for a tuple $B$, which differs from $A$ by one polytope. The following Lemma considers a special case, where such a reduction can be done.

\begin{lem}\label{tuple_change}
Let $A$ be an essential $n$-tuple in $\R^n$ and $B$ be the tuple resulting from $A$ by substitution of $A(i)$ for $A(j)$ for some $i,j$. If $\prod A = 1$, then $\prod B = 1$.
\end{lem}

\begin{prf}
Consider the tuple $C$ resulting from $A$ by substitution of $A(j)$ for $A(i)$. By the Aleksandrov--Fenchel inequality~(\ref{AF}), we have
\begin{equation}\label{AH-F}
1 = (\prod A)^2\geq \prod B\cdot \prod C.
\end{equation}
Since $A$ is essential, the tuples $B$ and $C$ are at least linearly independent. Hence, the factors of the right-hand side of (\ref{AH-F}) are strictly positive and integer, so $\prod B = \prod C = 1$.
\end{prf}

\begin{rmrk}
We will constantly use the following approach while reducing the proof of Theorem~\ref{Thm1} to simpler cases. In the settings of Lemma~\ref{tuple_change}, assume that we need to prove Theorem~\ref{Thm1} for the tuple $A$. Assume also that one of $B,C$ is essential, say, $B$. If Theorem~\ref{Thm1} is proved for $B$, it means that all the polytopes of $A$ except possibly $A(i)$ are contained, up to translation, in a volume 1 lattice simplex $\D$. Then we apply Lemma~\ref{step_lem} and conclude that the polytope $A(i)$ is also contained in $\D$ up to a translation, which completes the proof.
\end{rmrk}

The following statement helps to reduce the dimension $n$ in the proof of Theorem~\ref{Thm1}.

\begin{lem}\label{essentiality_funct}
Let $A$ be an essential $n$-tuple in $\R^n$ and $B\subset A$ be
a subtuple. Assume $B$ has dimension $1$. Then the tuple $C=B\cup \{X\}$ is also essential, where $X$ is the mixed fiber body of $A\setminus B$ in the affine span of $\sum B$.
\end{lem}

\begin{prf}
Assume $C$ is not essential: there exists a proper subtuple $C'\subset C$ such that $\dim C'\leq 0$. We have actually $X\in C'$ and $\dim C'=0$, otherwise the tuple $A$ would not be essential. Consider an arbitrary $Y\in C\setminus C'$. Let $I\subset \sum C'$ be a segment. Consider $A' = A\setminus \{Y\}$, $A_{\l} = A'\cup \{Y+ \l I\}$, $C_{\l} = (C\setminus \{Y\})\cup \{Y+ \l I\}$, $\l\geq 0$. We have
\begin{equation}\label{prod_lambda_1}
\prod A_{\l} = \prod C_{\l} = \prod C' \cdot \prod C'' = \prod C,
\end{equation}
where $C''$ is the projection of the tuple $C_{\l}\setminus C'$ along the affine span of
$\sum C'$, which does not depend on $\l$. On the other hand, monotonicity property and linearity imply respectively that
\begin{equation}\label{prod_lambda_2}
\prod A_{\l} \geq \prod (A'\cup \{\l I\}) = \l \prod (A'\cup \{I\}).
\end{equation}
Equations~\ref{prod_lambda_1},~\ref{prod_lambda_2} imply that $\prod (A'\cup \{I\}) = 0$, which contradicts to the statement of Lemma~\ref{essentiality_funct} requiring $A$ to be essential.
\end{prf}

\section{Proof of Theorem \ref{Thm1}, existence}\label{exists}

In this section we prove the main statement of Theorem~\ref{Thm1}. That is:

\begin{prop}
If the mixed volume of en essential $n$-tuple $A$ in $\R^n$ is 1, then all the polytopes $A_j$ of the tuple $A$ are contained, up to translations, in a volume 1 lattice simplex.
\end{prop}

We construct a multi-level induction on the following four aspects of the tuple $A$ that we list in the decreasing order of their magnitudes, from major aspects to minor ones.
\begin{enumerate}
\item The prime aspect is the dimension $n$ of the ambient space $\R^n$ containing the polytopes of $A$. The main goal of induction is to reduce the proof of Theorem~\ref{Thm1} to the cases of smaller dimensions. The extreme case of $n=1$ is obvious.
\item The secondary aspect is the number $k$ of distinct polytopes in $A$. Two polytopes $X_1,X_2\subset \R^n$ are considered as distinct, if $X_1\neq X_2$. We reduce the proof to smaller $k$. The extreme case $k=1$ means that all the polytopes $A_j$ coincide, and $\Vol(A_j)=\prod A = 1$. Thus $A_j$ is a volume 1 lattice simplex itself.
\item The third aspect denoted by $g=g(A)$ is rather technical and is defined as follows. If there exists $j$ such that $\dim (A(j)+A(j'))< \dim A(j) + \dim A(j')$ for some $j'<j$, then $g(A)$ is the minimum among those $j$. If there is no such $j$, then we assume $g(A)=1$. As we show below, the cases $g=1, g=2$ are technically more convenient, so we reduce the prove to smaller values of $g$.
\item At last, we consider the multiplicity vector $a = (a(1),a(2),\ldots,a(n))$ of the tuple $A$. Here $a(j)$ is the number of polytopes in $A$ that are equal to the polytope $A_j$ up to translations. For example, if $k=1$, then $a(j)=n$ for all $j=1,2,\ldots, n$. If $n=3$, $A_1=A_3$ and $A_2$ does not coincide with $A_1$ up to a translation, then $a(1)=a(3)=2$ and $a(2)=1$. A tuple $A$ is said to be richer than a tuple $B$ if the multiplicity vector of $A$ is smaller than that of $B$ in the sense of lexicographic order: $A\rhd B$. Respectively, $B$ is said to be poorer than~$A$: $B\lhd A$. As one can see below, we inductively reduce the proof of Theorem~\ref{Thm1} to poorer tuples.
\end{enumerate}

Base of induction is provided by the cases $n=1$ and $k=1$, both are vacuous, as we noted above.

We provide here a brief description of the way our induction is generally conducted. For a tuple $A$, suppose there exist $j$ such that the multiplicity $a(j)$ of $A_j$ is less than $\dim A_j-1$. We consider a tuple $B$ obtained from $A$ via substitution of some polytope of $A$ for a polytope equal to $A_j$ up to a translation. If $A\rhd B$, we have a chance to reduce the proof to the tuple $B$. However, there are two different cases.

First, assume $B$ is not essential. In this case, we show how to find a subtuple $A'\subset A$ of $\dim A'=1$ which contains at least two polytopes that do not coincide up to a translation. By induction on $n$, we manage to conclude that all the polytopes of $A'$ are contained, up to translations, in a volume~1 lattice simplex $\delta$ of dimension $d=|A'+1|$. Then we consider the tuple $C$ obtained from $A$ by substitution of all the polytopes of $A'$ for the simplex $\d$. Using Proposition~\ref{prop-2}, we conclude that $\prod C=1$. Then we reduce the proof of Theorem~\ref{Thm1} for $A$ to its proof for $C$, a tuple with smaller value of $k$.

Second, assume $B$ is essential. Then we reduce the proof for $A$ to the proof for $B$ by the arguments proposed in the remark below Lemma~\ref{tuple_change}.

This way, while we change some polytopes of $A$ to another ones, we tend to increase the dimensions of polytopes and decrease their variability in the tuple $A$. However, it turns out no longer possible at some step. Then we face the most difficult point of the proof. We analyze it below in the case~4 of induction step. The next subsection is devoted to a rigor description of the proof.

\subsection*{Step of induction}

In the sequence $A_1,A_2,\ldots, A_n$, there are exactly $k$ members $A_{i_1}, A_{i_2}, \ldots, A_{i_k}$ such that $A_{i_j}\neq A_{s}$ for all $s<i_j$. In other words, $1=i_1<i_2<\ldots< i_k$ identify the former representatives $A_{i_j}$ of~$k$ distinct polytopes in the tuple $A$. We also set $i_0=0$ for convenience.
Let $l$ be the maximal element of $\{0,1,\ldots,k-1\}$ such that $a(i)=\dim A_i-1$ for all $1\leq i\leq i_l$. It may be actually the case, where $a(i)=\dim A_i-1$ for all $i$ including $i_k$. In this case we assume $l=k-1$, since it is convenient for further reasoning. We split the induction step into several cases as follows.

1. We have $g(A)>2$. Consider arbitrary $g',g''\in \{1,2,\ldots, n\}$ such that $\dim (A_g+A_{g'})< \dim A_g+\dim A_{g'}$. Reorder the polytopes of $A$ in such a way that the first two values are $A_{g'}, A_{g''}$ and apply Theorem~\ref{Thm1} to the obtained tuple, which holds by induction.

2. There is an $n'$-subtuple $A'\subset A$ with $n'\leq n-2$ and $\dim A' = 1$ that contains at least two distinct polytopes. Let $\AA$ be the mixed fiber body of $A\setminus A'$ in the affine span $L$ of the polytope~
$\sum A'$. We have $\prod (A'\cup \AA) = \prod A = 1$. Lemma~\ref{essentiality_funct} imply also that the tuple $A'\cup \AA$ is essential. Since $\dim L =n'+1< n$, statement of Theorem~\ref{Thm1} holds for $A'\cup \AA$ by induction, thus there is a volume 1 lattice simplex $\d\subset L$ containing all the polytopes
of $A'$, as well as the polytope $\AA$, up to a translation.

Let $C'$ be the tuple formed by $|A'|$ copies of $\d$. Consider the tuple $C: = (A\setminus A') \cup C'$. Due to the monotonicity property, we have $1= \prod (A'\cup \AA) \leq \prod (C'\cup \AA) \leq \Vol(\d) = 1$, and thus $\prod C = \prod (C'\cup \AA) = 1$. As the tuple~$A$ is essential, the tuple $C$ is also essential. Note that the values of~$A\setminus A'$ do not coincide with those of $A'$. In fact, if there were a polytope $X\in A\setminus A'$ that coincide with one of the polytopes of $A'$, then dimension of subtuple  $(A'\cup X)\subset A$ would be $\dim A'-1 =0$, and thus tuple $A$ would not be essential. Therefore, we conclude that the number of distinct values of $C$ is less than $k$, since $C'$ has one value $\d$ for at least two distinct values in~$A'$. Consequently, the statement of Theorem~\ref{Thm1} holds for $C$ by induction and therefore holds also for $A$.

3. The cases 1, 2 do not occur, and we have $l<k-1$. Let $B$ be the $n$-tuple obtained from $A$ by substitution of the last representative of the polytope $A_{i_{k}}$ in the sequence $A_1,A_2,\ldots,A_n$ for the polytope $A_{i_{l+1}}$. Lemma~\ref{tuple_change} implies that $\prod B = 1$. We can also observe $g(B)\leq g(A)$ in both cases $g(A)=1$ and $g(A)=2$.

We claim that $B$ is essential. If not, there exists a proper subtuple $B'\subset B$ containing $a(i_{l+1})+1$ copies of $A_{i_{l+1}}$ with $\dim B'=0$. Consider the subtuple $A'$ of $A$ that is obtained from  $B'$ by removing one representative of $A_{i_{l+1}}$. We have $\dim A' = 1$. Since $a_{l+1}< \dim A_{i_{l+1}}-1$, the tuple $A'$ contains a polytope
different from $A_{i_{l+1}}$, and therefore the case 2 occur.

The number of distinct values of $B$ is less or equal $k$, $g(B)\leq g(A)\leq 2$, while $B\lhd A$, since the multiplicities $b(i)$ of the polytopes $B_i$ in $B$ coincide with $a(i)$ for $i<i_{l+1}$ and $b(i_{l+1}) = a(i_{l+1})+1$. Therefore, statement of Theorem~\ref{Thm1} holds for $B$ by induction: there is a volume 1 lattice simplex $\D\in \R^n$ containing all $B_i$ and thus also all the polytopes $A_i$, possibly except for $A_{i_{k}}$, up to translations. 
Applying Lemma~\ref{step_lem}, we conclude that $\D$ contains also 
$A_{i_{k}}$ up to a translation, which completes the proof.

4. The cases 1, 2 do not occur, and $l=k-1$. This one is the most difficult case, and we need two more lemmas to cope with it. The first one is a technical tool, and the second one covers the statement of Theorem~\ref{Thm1} for a distinguished special case.

\begin{dfntn}
The polytopes $A_1,A_2,\ldots, A_k$ of a tuple $A$ in $\R^n$ are said to be {\it transversal}, if  $\dim \sum A_i = \sum \dim A_i$. A polytope
$A$ is said to be {\it transversal to}
a set of polytopes $B_1,B_2,\ldots, B_k$, if the polytopes
$A,B_1,\ldots, B_k$ are transversal.
\end{dfntn}

\begin{lem}\label{diffucult_step}
Assume an $n$-tuple $A$ in $\R^n$ contains transversal volume 1 lattice simplices
$$
{\d_1,\d_2,\ldots, \d_p\subset \R^n},
$$
each $\d_j$ with multiplicity $a(j) = \dim \d_j -1$. Assume $A$ contains also a polytope~$X$, which is transversal to each of $\d_j$ but not transversal to the set $\d_1,\d_2,\ldots, \d_p$. Assume $\prod A = 1$, and $\prod (A'\cup \{\d_j\}) = 1$, $j=1,2,\ldots,p$, where $A':=A\setminus \{X\}$. Then, there exist $s_1,s_2$ and a volume~1 lattice
simplex $\D$ adjacent to $\d_{s_1},\d_{s_2}$ up to translations
such that $\dim \D - \dim \d_{s_1} - \dim \d_{s_2}$ equals 0 or 1, and $\prod (A'\cup \D)=1$.
\end{lem}

\begin{prf}
Let $\Ps\subset \ZZ$ be the set of covectors $\a$ such that $\prod (A')^{\a}>0$. We use the following
\begin{prop}\label{covectors}
For each $j=1,2,\ldots, p$:
\begin{enumerate}
\item There are covectors $\a^j_1,\ldots, \a^j_{l_j}\in \Ps$, where $l_j=\dim \d_j+1$, such that $\prod (A')^{\a^j_s} = 1$ and $\d_j^{\a^j_s}$, $s=1,2,\ldots, l_j$, are the $l_j$ facets of $\d_j$.

\item We have $\d_j^{\a} = \d_j$ for any covector $\a\in \Ps\setminus \{\a^j_1, \a^j_2, \ldots, \a^j_{l_j}\}$, .
\end{enumerate}
\end{prop}

\begin{prf}
Let $B_j\subset A$ be the subtuple consisting of $a(j)$ copies of
$\d_j$.  Consider the projection $C_j$ of the tuple $A'\setminus B_j$ along the affine span $L_j$ of the simplex $\d_j$. We have $\prod C_j = \prod C_j \cdot \Vol(\d_j)= \prod (A'\cup \{\d_j\}) = 1$. Let
$v_1,v_2,\ldots, v_{l_j}$ be the vertices of the simplex $\d_j$. For each edge $e_{st}=\{v_s,v_t\}, s\neq t,$ we have
\begin{equation}\label{prod_1}
\prod (A'\cup \{e_{st}\}) = \prod (B_j\cup \{e_{st}\})\cdot \prod C_j = 1.
\end{equation}
On the other hand, equation~(\ref{IV_reduce}) implies:
\begin{equation}\label{prod_2}
\prod (A'\cup \{e_{st}\}) = \sum_{\a\in \Psi} \left ( \max(\a(v_s),\a(v_t))\cdot \prod (A')^{\a} \right ).
\end{equation}
For a coordinate system in $\R^n$ centered in $v_t$, the summands of the right-hand side of~(\ref{prod_2}) are non-negative and integer.
Therefore, the equations~(\ref{prod_1}) and~(\ref{prod_2}) imply: there is a covector $\a_{st}\in \Ps$ with $\a_{st}(v_s)=1$, $\prod (A')^{\a_{st}} = 1$ and for any other covector $\a\in \Psi$, we have $\a(v_s)=0$. Observe that
\begin{equation}\label{prod_3}
1= \prod (A'\cup \{\d_j\}) = \sum_{\a\in \Psi} \left ( l_{\d_j}(\a)\cdot \prod (A')^{\a} \right ),
\end{equation}
thus, for the coordinate system centered in $v_t$, there is exactly one covector $\a^j_t\in\Ps$ such that $\max (\a^j_t)|_{\d_j})>0$.  In particular, $\a_{st} = \a^j_t$ not depending on $s$, and thus we have $\a^j_t(v_t)=0$ and $\a^j_t(v_i)=1$ for each $i\neq t$. Therefore, $\d_j^{\a^j_t}$ is the facet of $\d_j$ that do not contain $v_t$. Equation~\ref{prod_3} also implies that $\prod (A')^{\a^j_t}=1$ and this completes part 1. For part 2, we apply the same arguments using~(\ref{prod_3}).
\end{prf}

We continue the proof of Lemma~\ref{diffucult_step} with the two following cases.

1. Assume we have $\a^u_{s_1}=\a^v_{s_2}=\a_0$ for some $u\neq v$ and $s_1,s_2$. Assume for simplicity that the vertices of the polytopes $\d_u, \d_v$ not belonging to $\d_u^{\a^u_{s_1}}, \d_v^{\a^v_{s_2}}$ coincide with the origin (it is actually true up to translations). Consider $\D=\langle\d_u\cup \d_v \rangle$ to be the convex hull of $\d_u\cup \d_v$. In accordance with Proposition~\ref{covectors}, equation~\ref{IV_reduce} implies that $\prod (A'\cup \D) = \sum_{\a\in \Ps} l_\D(\a) \prod (A')^{\a} = l_\D(\a_0)\cdot 1 = 1$. In this case, we have $\dim \D - \dim \d_u - \dim \d_v = 0$.

2. Now assume the contrary: we have $\sum_{j=1}^p l_j$ different covectors $\a^j_s$. Assume for simplicity that each of $\d_j$, $j=1,2,\ldots, p$, contains the origin. Since $\dim (X+ \sum_{j=1}^p \d_j) < \dim X + \dim \sum_{j=1}^p \d_j$, there exists a non-degenerate segment $[x,y]\subset X$ such that $z=x-y\in \oplus_{i=1}^p L_j$, that is, $z = \sum z_j$, where $z_j$ belongs to the affine span $L_j$ of $\d_j$. There are at least two values $u,v$ of the variable $j$ such that $z_j\neq0$, otherwise, the polytope $X$ would not be transversal to those $\d_j$ corresponding to the unique $z_j\neq 0$. Due to Proposition~\ref{covectors}, part~2, we have $\a^{j_1}_s(z_{j_2})=0$ for any $j_1\neq j_2$, and thus $\a^j_s(z) = \a^j_s(z_j)$. There exist $s_1,s_2$ such that $\a^u_{s_1}(z_u)>0$ and $\a^v_{s_2}(z_v)>0$. This means that $\a^u_{s_1}(x)>\a^u_{s_1}(y)$, $\a^v_{s_2}(x)>\a^v_{s_2}(y)$, and thus the covectors $\a_1=\a^u_{s_1}, \a_2=\a^v_{s_2}$ are not constant on the polytope $X$. For each of the two covectors, the range of values on $X$ is a unit segment because of the following equation:
\begin{equation}\label{prod_X}
\sum_{\a\in \Ps} \left ( l_X(\a)\cdot \prod (A')^{\a} \right ) = \prod A.
\end{equation}

First, assume there is a vertex $x\in X$ with $\a_1(x)=\min(\a_1|_{X})$, $\a_2(x)=\min(\a_2|_{X})$. W.l.o.g., assume that $x$ coincides with the origin. The left-hand side summands of equation~(\ref{prod_X}) are non-negative, while those for $\a=\a_1, \a_2$ are positive and integer.
Therefore, $\prod A \geq 2$, which contradicts to the conditions of Lemma~\ref{diffucult_step}.

Thus it follows such a vertex $x$ does not exist. In this case, there exist two different vertices $x_1,x_2\in X$ such that $\a_1(x_2)=\a_1(x_1)+1$, $\a_2(x_1)=\a_2(x_2)+1$. Denote $x=x_2-x_1$. We claim that $\a(x)=0$ for $\a\in \Psi, \a\neq \a_1, \a_2$. If not, there exist $\a_0\in \Psi, \a_0\neq \a_1,\a_2$ with $\a_0(x_2)>\a_0(x_1)$ (or $\a_0(x_1)>\a_0(x_2)$). Translate
the polytope $X$ in such a way that $x_1$ coincides with the origin (or $x_2$ coincides with the origin). The left-hand side summands of equation~(\ref{prod_X}) become non-negative, while those corresponding to $\a=\a_1,\a_0$ (or $\a=\a_2,\a_0$) are strictly positive. Therefore, $\prod A \geq 2$, which is a contradiction.

Translate the polytopes $\d_u, \d_v$ in such a way that their vertices not belonging to $\d_u^{\a_1}, \d_v^{\a_2}$ coincide with the origin. Next, translate the polytope
$\d_{v}$ by adding vector $x$. Consider the simplex $\D=\langle \d_u\cup \d_v \rangle$. We have $\max(\a^u_s|_\D)=0$ for $s\neq s_1$, $\max(\a^v_s|_\D)=0$ for $s\neq s_2$, $\max(\a^u_{s_1}|_\D)=1$, $\max(\a^v_{s_2}|_\D)=0$, and $\max(\a|_\D)=0$ for other covectors $\a\in \Psi$. Therefore we have $\prod (A'\cup \D)= \sum_{\a\in \Psi} \left (l_\D(\a)\cdot \prod (A')^{\a} \right ) =1$. In this case, $\dim \D - \dim \d_u - \dim \d_v = 1$.
\end{prf}

\begin{lem}\label{integral} Let $X\subset \R^n$ be a lattice polytope, and $\d\subset \R^n$ be a volume 1 lattice
simplex, $\dim \d=d<n$. Let
$A$ be the tuple consisting of $(d-1)$ copies of $\d$ and $n-d+1$ copies of $X$. Assume $B=(A\setminus \{X\})\cup \{\d\}$. If $\prod A = \prod B = 1$, there exists a volume 1 lattice
simplex $\D$ that contains the polytopes
$\d, X$ up to translations.
\end{lem}

\begin{prf}
Proposition~\ref{prop} implies that the tuple $B$ is linearly independent and thus the projection $Y$ of the polytope
$X$ along the affine span of the simplex $\d$ has dimension $n-d$. Therefore, there exists a volume 1 lattice
simplex $\d'\subset X$ transversal to
$\d$ such that $\dim \d' = n-d$. Since $A$ is linearly independent, we have $\dim X\geq n-d+1$, and therefore $X$ is not contained in the affine span of $\d'$. Thus there exists a lattice segment $X'\subset X$, which is transversal to $\d$, as well as to $\d'$. Consider the tuple $C$ consisting of $d-1$ copies of
$\d$, $n-d-1$ copies of
$\d'$, and one copy of each of $X$ and $X'$. Assume
$C'=C\setminus \{X'\}$, $C_1=C'\cup\{\d'\}$, $C_2=C'\cup\{\d\}$. The tuples $C, C_1, C_2$ are linearly independent.
Since the monotonicity property implies that $\prod C_1, \prod C \leq \prod A = 1$, $\prod C_2\leq \prod B = 1$, we actually observe $\prod C = \prod C_1= \prod C_2 =1$. Applying Lemma~\ref{diffucult_step} to the tuple~$C$, we obtain a volume 1 lattice full-dimensional
simplex $\D$ in $\R^n$ such that $\prod (C'\cup \{\D\})=1$ and the simplices $\d,\d'$ are contained in $\D$ up to translations. Note that the tuple $C'\cup \{\D\}$ is essential, and thus Lemma~\ref{step_lem} implies that the polytope $X\in C'$ is also contained in $\D$ up to a translation.
\end{prf}

Now we return to the case 4 of induction step. Assume $\dim A_{i_{j_0}} = n$ for some $j_0\in \{1,2\ldots, l\}$. Then $a(i_{j_0})=n-1$, and we have actually $k=2$. Lemma~\ref{tuple_change} implies that the volume of $A_{i_{j_0}}$ equals 1, and thus the polytope
$A_{i_{j_0}}$ is a volume 1 lattice
simplex $\D$. Lemma~\ref{step_lem} implies
also that the only polytope
of $A$ different from $A_{i_{j_0}}$ is also contained in $\D$ up to a translation, which completes the proof.

In what follows, we assume that $\dim A_{i_j} < n$ for $j=1,2,\ldots,l$. Let us show that each of the polytopes
$A_{i_j},\, j=1,2\ldots, l,$ is a simplex. Since the affine span $L_j$ of the polytope $A_{i_j}$ has dimension less than $n$, Theorem~\ref{Thm1} holds there by induction. Let $B_j$ be the tuple consisting of $a(i_j)$ copies of $A_{i_j}$. Let $\AA_j\subset L_j$ be the mixed fiber body of the set of the
other polytopes of $A$. Applying Theorem~\ref{Thm1} to the tuple $B_j\cup \{\AA_j\}$, we get a volume 1 lattice simplex  $\d_j\in L_j$, which contains $A_{i_j}$. Since $\dim A_{i_j} = \dim \d_j,$ we get $A_{i_j} = \d_j$.

Note that $l=k-1>0$, and, therefore, $A_{i_1}$ is a simplex. We consider the two following subcases below:

A). We have $g=2$. If $k=2$, we apply Lemma~\ref{integral}. Assume $k>2$. In this case, $l\geq 2$, and thus there is a subtuple $C\subset A$ consisting of $a(i_1)$ copies of $A_{i_1}$ and $a(i_2)$ copies of $A_{i_2}$. Since $g=2$, we have $i_1=1$, $i_2=2$, and $\dim (A_{i_1}+A_{i_2}) < \dim A_{i_1} + \dim A_{i_2}= a(i_1)+a(i_2)+2$. Therefore, $\dim C< 2$ and, since $A$ is essential, we have $\dim C = 1$. If we had also $a(i_1)+a(i_2) < n-1$, then the case 2 would take place. Thereby, $a(i_1)+a(i_2) = n-1$, and $k=3$. We apply Lemma~\ref{integral} to the tuple $B:=(A\setminus\{A_{i_3}\})\cup \{A_{i_2}\}$ and obtain a  volume~1 lattice simplex $\D$, which contains the polytopes of $B$ up to translations. Since the tuple $A$ differs from $B$ in a single polytope, Lemma~\ref{step_lem} completes the proof.

B). We have $g=1$. Let $p$ be the maximal element of \{1,2,\ldots, k\} such that the polytopes $A_{i_1},A_{i_2},\ldots, A_{i_p}$ are transversal. We claim that $p\leq l=k-1$. In fact,
$$
n = \sum_{j=1}^k a(i_{j}) < \sum_{j=1}^k \dim A_{i_j}.
$$
Taking into account that $n = \dim \sum_{j=1}^k A_{i_j}$, we get $\dim \sum_{j=1}^k A_{i_j} < \sum_{j=1}^k \dim A_{i_j}$, thus $p\neq k$.

Consider the polytope $X = A_{i_{p+1}}$ and the $(n-1)$-tuple $A':=A\setminus \{X\}$. Since $g=1$, the polytope $X$ is transversal to each of $A_{i_j}, j\leq p$. Lemma~\ref{diffucult_step} provides a  volume 1 lattice simplex $\D$ such that $\prod (A'\cup \{\D\})=1$, $A_{i_u},A_{i_v}$ are contained in $\D$ up to translations for some $u,v\leq p$ and $\dim \D - \dim A_{i_{u}} - \dim A_{i_{v}}=d$, which equals 0 or 1. We claim that $B=A'\cup\{\D\}$ is essential. Otherwise there is a subtuple $C\subset B$ such that $A_{i_u}, A_{i_v}, \Delta\in C$, $\dim C=0$, and $|C|\leq n-1$. Then $C'=C\setminus \{\Delta\}$ is a subtuple of $A$ containing at least two distinct polytopes $A_{i_u}, A_{i_v}$, and we have $\dim C'=1$, $|C'|\leq n-2$, which implies that case~2 takes place.

Let $B'\subset B$ be the subtuple consisting of $a(i_u)$ copies of
$A_{i_u}$, $a(i_v)$ copies of $A_{i_v}$ and the polytope $\D$. The polytopes of $B'$ are contained in the volume 1 lattice simplex $\D$, and we have $\dim B' = 1$. Consider the mixed fiber body $\widetilde{B}$ of $B\setminus B'$ in the affine span of $\D$. Lemma~\ref{essentiality_funct} implies that the tuple $B''= B'\cup\{\widetilde{B}\}$ is essential. We have also $\prod B'' = \prod B=1$, thus Lemma~\ref{step_lem} implies that $\widetilde{B}$ is contained in $\D$ up to a translation. Let $C$ be the tuple formed by $|B'|$ copies of $\D$. The tuple $E=(B\setminus B')\cup C$ is essential, since the tuple $B$ is. We have also $\prod E = \prod C\cup\{\widetilde{B}\}=1$. The number of distinct values of $E$ is lesser than $k$. Therefore, Theorem~\ref{Thm1} holds for $E$ by induction and thus holds for~$B$. Since $A$ differs from $B$ in a single polytope, Lemma~\ref{step_lem} completes the proof.

\section{Proof of Theorem \ref{Thm1}, uniqueness} \label{unique}

\begin{lem}\label{pullback1}
Let $B$ be an essential tuple of $m$ subsimplices in the standard simplex in $\Z^m$, and assume it is the projection of a tuple $A$ in $\Z^m\oplus\Z^1$ such that $\dim\sum A=m+1$. There exists at most one tuple of numbers $c_2,c_3,\ldots,c_m$ such that $A(1),A(2)+(0,\ldots,0,c_2),\ldots,A(m)+(0,\ldots,0,c_m)$ are contained in a volume 1 simplex.\end{lem}


\begin{rmrk}  If we know that such $c_2,c_3,\ldots,c_m$ do exist, then they can be found in polynomial-time as follows. Denote the vertices of the standard simplex in $\Z^m$ by $v_0,v_1,\ldots,v_m$. For every $k=0,1,\ldots,m$, verify if one can translate the elements of $A$ into a volume 1 simplex, whose edge is contained in the line $\{v_k\}\times\Z^1$. This can be verified in polynomial-time as follows. Reorder $A$ and $B$ so that, for every $i=2,3,\ldots,m$, we have $B(i)\cap \Bigl(B(i-1)\cup\ldots\cup B(1)\Bigr)\setminus\{v_k\}\ne\varnothing$, then, for every $i=2,3,\ldots,m$, translate $A(i)$ along $\{0\}\times\Z^1$ so that $A(i)\cap \Bigl(A(i-1)\cup\ldots\cup A(1)\Bigr)\setminus(\{v_k\}\times\Z^1)\ne\varnothing$, then verify if the vertices of the resulting translated copies of $A(1),A(2)\ldots,A(m)$ are the vertices of a volume 1 simplex. \end{rmrk}

{\it Proof of lemma.} Assume (without loss in generality) that there exists one such collection $c_2=\ldots=c_m=0$, the resulting simplex in $\Z^m\oplus\Z^1$ is the standard one, and $A(1)\ni 0$. We now try to find another collection $(c_2,\ldots,c_m)\ne 0$, such that $A(1),A(2)+(0,\ldots,0,c_2),\ldots,A(n)+(0,\ldots,0,c_n)$ is contained in another volume 1 simplex.

If $(0,\ldots,0,1)\in A(1)$, this is obviously impossible.
Otherwise we can assume without loss in generality that $(0,\ldots,0,1)\in A(m)$, thus $c_m=-1$.
Since $B$ is essential, the union $B(1)\cup\ldots\cup B(m-1)$
is connected, and its convex hull is the standard simplex in $\Z^m$. This implies that $c_2=\ldots=c_{m-1}=0$.
On the other hand, since $B$ is essential, $B(m)$ contains at least three vertices $0, a, b$ in $\Z^m$. This imples that the union of $A(1),A(2)+(0,\ldots,0,c_2),\ldots,A(m)+(0,\ldots,0,c_m)$ contains at least $m+3$ points $(0,0,\ldots,0),(1,0,\ldots,0),\ldots,(0,\ldots,0,1,0),(a,-1),(b,-1)$, which are not contained in a volume 1 simplex in $\Z^m\oplus\Z^1$. $\quad\Box$

We say that a tuple $A$ is $k$-unique, if, for some $i_1,\ldots,i_k$,
there exists a unique, up to a simultaneous translation, tuple of translations of the polytopes $A(i_1),\ldots,A(i_k)$ that can be extended to a tuple of translations of $A(1),\ldots,A(n)$ into a lattice simplex of volume 1. 
In particular, any essential tuple of mixed volume 1 is $1$-unique.
The rest of this section is devoted to the proof of the following statement, equivalent to the uniqueness in Theorem 2:

\begin{prop} Assume that an essential $n$-tuple $A$ in $\Z^n$ has mixed volume 1. If $A$ is $k$-unique, $k\geqslant 1$, then $A$ is $(k+1)$-unique.\end{prop}

{\it Proof.} With no loss in generality, we can assume that $i_1=1,\ldots,i_k=k$. Moreover, replacing each of $A(1),\ldots,A(k)$ with the convex hull of the union of their translated copies in a volume 1 simplex, we can also assume that $A(1)=\ldots=A(k)$.

For an edge $e\subset A(1)$, consider an epimorphism $\pi:\Z^n\to\Z^{n-1}$, such that $\pi(e)$ is a point. Denote $\pi (A(i))$ by $B(i)$ for $i=2,\ldots,n$, then the tuple $B$ has mixed volume 1. Let $v_{k+1},v_{k+2},\ldots,v_n$ be any set of vectors such that $A(1)=\ldots=A(k),A(k+1)+v_{k+1},\ldots,A(n)+v_n$ are contained in a volume 1 simplex. Then  $B(2)=\ldots=B(k),B(k+1)+\pi (v_{k+1}),\ldots,B(n)+\pi (v_n)$ are also contained in a volume 1 simplex, and one of the following two cases takes place.

1) $B$ is essential. Then Theorem~\ref{Thm1} holds for $B$ by induction on $n$, and thus
$\pi (v_{k+1}),\pi (v_{k+2}),$ $\ldots,$ $\pi (v_n)$
are uniquely determined by the tuple $A$ (up to a simultaneous addition of some vector to each of them). Then, applying Lemma~\ref{pullback1} to the tuple $B(2)=\ldots=B(k),B(k+1)+\pi (v_{k+1}),\ldots,B(n)+\pi (v_n)$, the vectors $v_{k+1},v_{k+2},\ldots, v_n$ are also uniquely determined by  $A$ (up to adding some vector $v$ to all of them), and 
the tuple $A$ is $n$-unique.

2) $B$ contains a non-empty essential subtuple $B'$ that consists of $B(i),\, i\in I\subset \{k+1,k+2,\ldots, n\}$, and satisfies $\dim B'=0$. Then we have $\prod B'=1$, and Theorem~\ref{Thm1} holds for $B'$ by induction on $n$. This implies that $\pi (v_i), i\in I$, are uniquely determined by $A$ (up to a simultaneous addition of some vector to each of them). Then, by Lemma~\ref{pullback1}, the vectors $v_i, i\in I$, are also uniquely determined by $A$, up to adding some vector $v$ to all of them. This vector $v$ is uniquely determined by the condition that the convex hull $S$ of the union of $A(i)+v_i,\, i\in I$, contains the edge
$e$ of $A(1)$ 
(note that $S$ has to contain $e$ up to translation, otherwise the restriction of $\pi$ to $S$ is injective, thus $|I|=\sum_{i\in I} A(i)$, and $A$ is not essential). Thus
the tuple $A$ is $(k+|I|)$-unique. $\quad\Box$

\section{Shifting polytopes into a unit simplex in polynomial time} \label{Sshift}

Algorithm (I-II) in Section 1 contains the following two steps, whose implementation in po\-ly\-no\-mial-time is not obvious.

1) {\it Given a tuple of Newton polytopes of positive mixed volume, find its maximal essential subtuple.} For any tuple $A$ of polytopes, the function $\dim I=|A|-\dim\sum_i A_i$ is submodular (see \cite{Sturmfels}), so a minimal subset $I$, satisfying $\dim I=0$ (i.e. corresponding to a maximal essential subtuple), can be found in polynomial-time (see, e.g., \cite{subm1}, \cite{subm2}). Note that, applying this algorithm to a tuple $A$ of mixed volume 0, we shall find $I$ such that $\dim I>0$, i.e. detect that the mixed volume of $A$ is not positive.

2) {\it Given an essential tuple of Newton polytopes $A$ of mixed volume 1, find how to translate the polytopes into a single volume 1 simplex.} Note that the first approach that comes to mind (choose a segment in every Newton polytope and try to assemble translated copies of these segments into a volume 1 simplex, then check if it contains the original Newton polytopes up to a translation) leads to an exhaustive search among $\sim$$n^n$ options, see Proposition \ref{cayley}. 
However, the following polynomial-time algorithm is provided by the reasoning in Section \ref{unique}.
We rely upon the notation introduced in the course of Section \ref{unique}.

\begin{enumerate}
\item Choose an arbitrary edge $e\in A(1)$ and consider the tuple $B$ that consists of the images of $A(i),\, i>1$, under the projection $\pi$ along $e$. For every $i_1<\ldots<i_p$, denote the tuple $B(i_1),\ldots,B(i_p)$ by $B\{i_1,\ldots,i_p\}$. Let $BI_1,\ldots,BI_s$ ($s$ may equal 0) be the essential subtuples in $B$, such that $\dim BI_j=0$, then it is crucial to notice that any $I$, such that $\dim BI=0$, contains one of $I_1,\ldots,I_s$, and these $s$ subtuples do not overlap, because $\dim B(I_i\cap I_j)=0$ for any $i$ and $j$.  
In particular, we have $|I_1|+\ldots+|I_s|<n$.

\item Applying the algorithm recursively to $BI_j$, we find vectors $u_i,\, i\in I_j$, such that the simplices $B(i)+u_i$ are contained in the same volume 1 simplex. Then, applying the remark after Lemma \ref{pullback1} to the simplices $B(i)+u_i,\, i\in I_j$, we find vectors $v_i,\, i\in I_j$, such that the simplices $A(i)+v_i$ and the edge $e$ are contained in the same volume 1 simplex $S_j$. Replacing every $B(i),\, i\in I_1\sqcup\ldots\sqcup I_s$, in the tuple $B$ with the convex hull of the union $\pi\bigl(A(1)\cup S_1\cup\ldots\cup S_s\bigr)$, we notice that the resulting tuple $B'$ is essential.

\item Applying the algorithm recursively to $B'$, we find vectors $u'_i,\,i>1$, such that the simplices $B'(i)+u'_i$ are contained in the same volume 1 simplex. Then, applying the remark after Lemma \ref{pullback1}, we find vectors $v'_i,\, i=1,\ldots,n$, satisfying $v'_i=v_i$ for $i\in I_1\sqcup\ldots\sqcup I_s$ and $\pi(v'_i)=u'_i$ otherwise, such that all the simplices $A(i)+v'_i$ are contained in the same volume 1 simplex.

\end{enumerate}
This algorithm (1-3) contains polynomially many operations in the
dimension $n$ and $s+1$ recursions to the dimensions
$|I_1|,\ldots,|I_s|,n-1$, satisfying $|I_1|+\ldots+|I_s|<n$. Thus it
runs in polynomial time.

Note that, applying this algorithm to an essential tuple $A$ of mixed volume $>1$, the sub-algorithm from the remark after Lemma \ref{pullback1} will fail (i.e. will give a non-simplex) at some step, so we shall detect that the mixed volume of $A$ is not 1.

\begin{exmpl} Let $e_1,e_2,e_3$ be the standard basis in $\R^3$, $a=a_1e_1+a_2e_2+a_3e_3\in\R^3$, and let $A$ consist of $(0,e_1,e_2),\, (0,e_2,e_3),\, (a,a+e_3,a+e_1)$. Then the algorithm proceeds as follows: $e=(0,e_1);\; B$ consists of $(0,e_2,e_3)$ and $(a_2e_2+a_3e_3, a_2e_2+a_3e_3+e_3)$; $s=1;\; BI_1$ consists of the latter segment; $u_3=-a_2e_2-a_3e_3$, $v_3=-a$; $B'$ consists of two copies of $(0,e_2,e_3)$, and $u_1=v_1=u_2=v_2=0$.
\end{exmpl}

Of course, there are lots of possibilities to optimize this algorithm.
For instance, in the case of dimension $<6$, one can always find $i$ and $j$ such that $A(i)$ and $A(j)$ have a common edge up to a translation. Thus, upon this translation, one can replace $A(i)$ and $A(j)$ with the convex hull of the union of them, and apply the same fact to the resulting tuple. However, this fails in higher dimensions:

\begin{exmpl} Let $0,\ldots,6$ be the vertices of the standard 6-dimensional simplex, then the mixed volume of the essential tuple of the triangles $012, 034, 056, 136, 145, 246$ equals 1, but no two of these triangles have a common edge. \end{exmpl}

\section{Further classification} \label{Smv2}

The same ideas as above allow to classify collections of lattice polytopes of larger mixed volume: in a forthcoming paper, we classify the polytopes, whose mixed volume equals two.
In the same way as for mixed volume 1, it is enough to classify the essential tuples.
Here is, for instance, the classification of essential tuples of mixed volume 2 in two and three dimensions. 

\begin{exmpl}\label{dim3}
{\rm Every essential pair of lattice polygons of mixed area 2 is contained in a pair that equals one of the following three (up to translations and an authomorphism of $\Z^2$):

\begin{center}
\noindent\includegraphics[width=10cm]{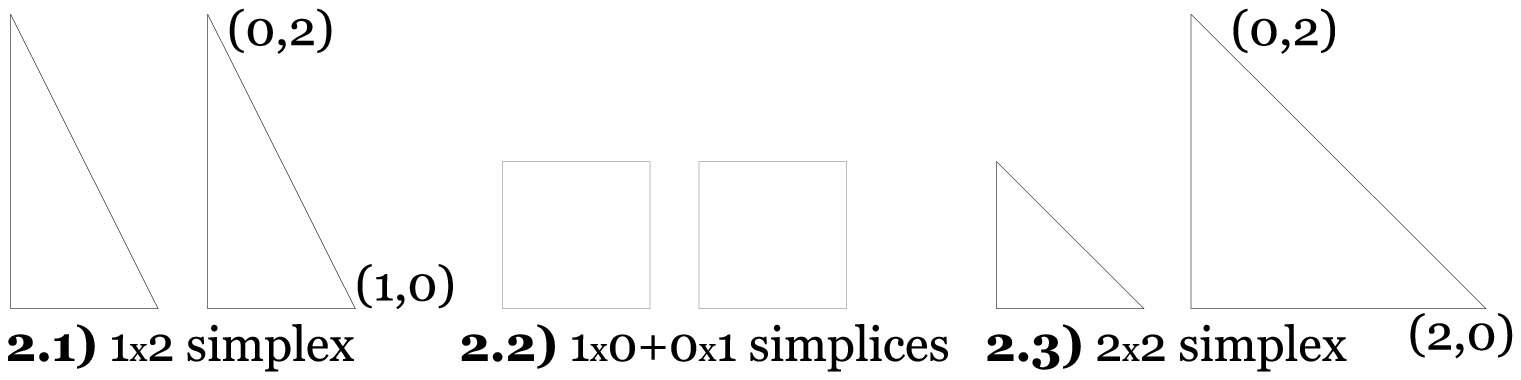}
\end{center}

Every essential triplet of lattice polytopes of mixed volume 2 is contained in a triplet that equals one of the following seven (up to translations and an authomorphism of $\Z^3$):

\begin{center}
\noindent\includegraphics[width=11cm]{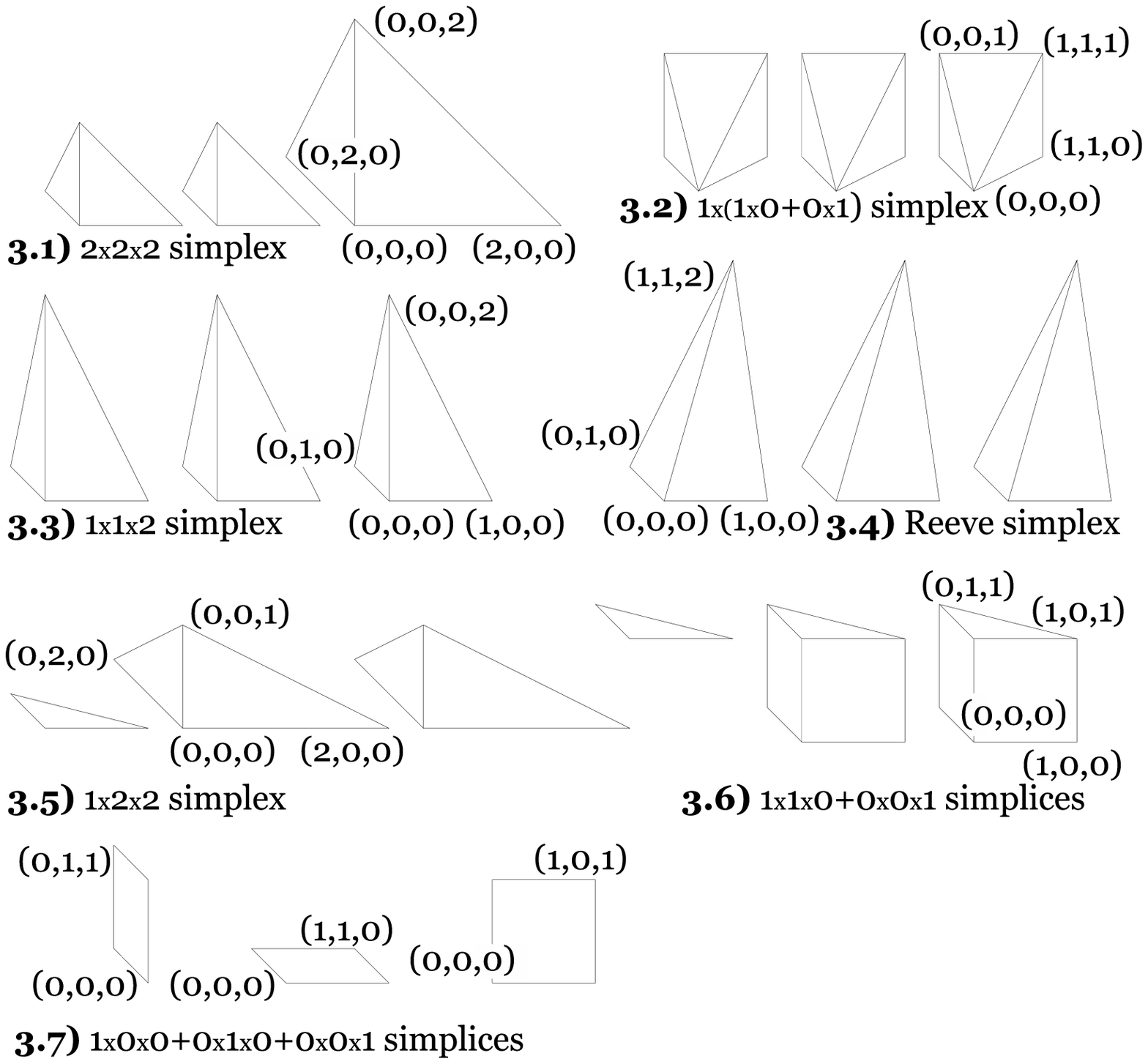}
\end{center}

Note that any system of polynomial equations in 2 or 3 variables, whose Newton polytopes constitute one of the tuples from the list above, can be solved by an explicit formula in quadratic radicals. The existence of such formula follows from general considerations, but producing it for particular multidimensional Newton polytopes is not always obvious, c. f. {\bf (3.7)} below. In dimension 3, the explicit formulas are as follows:

If a system $f_1(x,y,z)=f_2(x,y,z)=f_3(x,y,z)=0$ has the Newton polytopes as in {\bf (3.1)}, then, substituting $x$ and $y$ in $f_3=0$ in terms of $z$ from $f_1=f_2=0$, we obtain a quadratic equation on $z$.

If a system $f_1(x,y,z)=f_2(x,y,z)=f_3(x,y,z)=0$ has the Newton polytopes as in {\bf (3.2)},  {\bf (3.3)} or {\bf (3.4)}, then, adding suitable multiples of $f_3$ to $f_1$ and $f_2$, we obtain an equivalent system with the Newton polytopes as in {\bf (3.1)}.

If a system $f_1(x,y,z)=f_2(x,y,z)=f_3(x,y,z)=0$ has the Newton polytopes as in {\bf (3.5)} or {\bf (3.6)}, then, substituting $x$ in $f_2=f_3=0$ from $f_1=0$, we obtain a system in two variables with two solutions.

Finally, if a system $f_1(x,y,z)=f_2(x,y,z)=f_3(x,y,z)=0$ has the Newton polytopes as in {\bf (3.7)}, then the equations can be rewritten as $x=\varphi(y),\, y=\psi(z),\, z=\chi(x)$, where $\varphi,\, \psi$ and $\chi$ are M\"{o}bius transformations, so the problem is to find the fixed points of the M\"{o}bius transformation $x=\varphi(\psi(\chi(x)))$}.
\end{exmpl}

{\noindent(A. Esterov) National Research University Higher School of Economics. \newline Faculty of Mathematics NRU HSE, 7 Vavilova 117312 Moscow, Russia.}

{\noindent(G. Gusev) Moscow Institute of Physics and Technology (State University). \newline Department of Innovations and High Technology, 9 Institutskii per. \newline 141700 Dolgoprudny, Moscow Region, Russia.}

\end{document}